\documentclass[12pt,centertags,oneside]{amsart}

\usepackage{amsmath,amstext,amsthm,amscd,typearea,hyperref,stmaryrd}
\usepackage{amssymb}
\usepackage{a4wide}
\usepackage[mathscr]{eucal}
\usepackage{mathrsfs}
\usepackage{typearea}
\usepackage{charter}
\usepackage{pdfsync}
\usepackage[a4paper,width=16.2cm,top=3cm,bottom=3cm]{geometry}

\numberwithin{equation}{section}

\newtheorem{theorem}{Theorem}[section]
\newtheorem{definition}[theorem]{Definition}
\newtheorem{proposition}[theorem]{Proposition}

\newtheorem{lemma}[theorem]{Lemma}
\newtheorem{remark}[theorem]{Remark}

\newtheorem{problem}{Problem}

\newcommand{\cali}[1]{\mathscr{#1}}

\newcommand{\supp}{{\rm supp}}

\newcommand{\loc}{{loc}}
\newcommand{\ddc}{{dd^c}}

\newcommand{\ddbar}{{\partial\overline\partial}}

\newcommand{\ind}{{\bf 1}}

\newcommand{\Cc}{\cali{C}}

\newcommand{\Ec}{\cali{E}}

\newcommand{\Kc}{\cali{K}}

\renewcommand{\Mc}{\cali{M}}

\newcommand{\Rc}{\cali{R}}

\newcommand{\FS}{{\rm FS}}

\newcommand{\C}{\mathbb{C}}
\newcommand{\D}{\mathbb{D}}

\newcommand{\N}{\mathbb{N}}

\renewcommand{\P}{\mathbb{P}}

\newcommand{\TP}{{{\rm T}\mathbb{P}}}

\newcommand{\PGL}{{{\rm PGL}}}

\title{The Mandelbrot set is the shadow of a Julia set}

\author{Fran\c cois Berteloot}
\address{Universit\'e Toulouse 3, 
Institut Math\'ematique de Toulouse, 
118 route de Narbonne,
F-31062 Toulouse Cedex 9, France. }
\email{berteloo@picard.ups-tlse.fr}

\author{Tien-Cuong Dinh}
\address{Department of Mathematics, National University 
of Singapore, 10 Lower Kent Ridge Road, Singapore 119076. 
{\tt http://www.math.nus.edu.sg/$\sim$matdtc} }
\email{matdtc@nus.edu.sg} 

\date{}

\begin{document}

\maketitle

\begin{abstract}
Working within the polynomial quadratic family, we introduce a new point of view on bifurcations which naturally allows to see
the seat of bifurcations as the projection of a Julia set of a complex dynamical system in dimension three. 
We expect our approach to be extendable to other holomorphic families of dynamical systems.
\end{abstract}

\medskip\medskip


 \medskip


\tableofcontents

\section{Introduction} \label{s:intro}

Let  $f_c: M\to M$ be a family of holomorphic dynamical systems  which 
depend on a parameter $c\in \Sigma$ where both the phase space $M$ and the parameter space $\Sigma$
are complex manifolds. The bifurcation theory studies how the dynamics of $p_c:M\to M$ depends on the parameter
$c$. More specifically, we consider a subsystem $p_c:J_c\to J_c$
where $J_c$ is the support of a certain canonical ergodic invariant measure $\mu_c$ and can be seen as a \emph{Julia set}. We are interested in the set of parameters for which the dynamics of  $p_c:J_c\to J_c$ drastically changes under small perturbations of $c$ and call it the \emph{bifurcation locus}. The family of  ergodic dynamical systems $(p_c,J_c,\mu_c)_{c\in \Sigma}$ and its bifurcation locus are the central objects in the theory.

It is natural to focus on families which depend holomorphically on the parameter $c$ and, in that case, to consider
the holomorphic dynamical system $F:\Sigma\times M\to \Sigma\times M$ given by 
$F(c,z):=(c,p_c(z))$ instead of the  family $(p_c)_{c\in \Sigma}$. In this way, the parameters are included
in the phase space of the later system. This point of view 
has been followed by several authors to study bifurcations in projective holomorphic dynamics.
We refer to the lecture notes \cite{BCime, Ban, DSCime, dujardin} for these aspects.

The novelty in the present paper is to consider other dynamical systems induced by $F$ such as the associated dynamical system on the tangent bundle of $\Sigma\times M$. We expect  the bifurcation set, or at least a large set of bifurcation parameters, to naturally appear as Julia-type sets of those dynamical systems. 

We will restrict ourselves to the family of quadratic polynomials $p_c(z):=z^2+c$, which is both the first non-trivial 
and the most studied family. 
Before  precisely stating our result, let us introduce the main objects of the study.
We shall consider the projectivization $X:=\P\TP^2$ of the tangent bundle $\TP^2$ of $\P^2$ and denote by 
$\Pi$ the canonical projection 
$\Pi:X\to \P^2.$ The holomorphic map $F: \C^2\to \C^2$ given by $F(c,z)=(c,z^2+c)$, as well as its iterates $F^n$,
can be considered as rational maps from $\P^2$ to itself. We then lift them in a natural way to rational maps $\widehat F$ and $\widehat F^n$ from $X$ to $X$. These maps are induced respectively by the differentials of $F$ and $F^n$ which are rational self-maps of the tangent bundle $\TP^2$. 

For any $c\in \C$, the quadratic polynomial $p_c(z):=z^2+c$ is a holomorphic self-map of $\P^1$
whose  \emph{filled-in Julia set} and  \emph{Julia set}  are respectively denoted by $K_c$ and 
$J_c$. We recall that $J_c$ is equal to the boundary $bK_c$ of $K_c$ and $K_c$ is the set of points with bounded orbits in $\C$. 
We also have $K_c=\{z\in \C\;\colon \; g_c(z)=0\}$ where $g_c$ is the {\it dynamical Green function} of $p_c$ defined by
\begin{equation} \label{e:Green}
g_c(z):=\lim_{n\to\infty} 2^{-n} \ln^+ \vert p_c^n(z)\vert \quad \text{with}\quad p_c^n:=p_c\circ\cdots\circ p_c \text{ ($n$ times).}
\end{equation}
The {\it Mandelbrot set} $\Mc$ is defined by 
$$\Mc:=\{c\in \C\;\colon\; g_c(0)\le 0\}.$$
It is a classical result that the boundary $b\Mc$ of $\Mc$ is precisely the bifurcation locus of the quadratic polynomial family $(p_c)_{c\in \C}$. 

Define the two probability measures on the parameter space $\C$ and the phase space $\C$ of $p_c$ by 
$$m:=\ddc g_c(c) \quad \text{and} \quad \mu_c:=\ddc g_c,$$
where $\ddc:={i\over \pi} \ddbar$. It is well-known that $m$ and $\mu_c$ are respectively the equilibrium measures of $\Mc$ and $K_c$. We refer to \cite{CG,Si1} for properties of the Julia and Mandelbrot sets. We will consider $p_c$ as a dynamical system on the vertical line $\{c\}\times\C$ of $\C^2$. So both $J_c$ and $K_c$ are compact subsets of $\{c\}\times\C$ and hence of $\P^2$. The measure $\mu_c$ will also be identified with a probability measure on $\{c\}\times\C$ and can be seen as a probability measure of $\P^2$. So we can define a probability measure $\mu$ on $\P^2$ and a vertical  closed positive $(2,2)$-current $\Rc$ on $X$ by setting
$$\mu:=  \int  \mu_c\;dm(c) \quad \text{and} \quad \Rc:=\Pi^*(\mu).$$

The last current somehow provides a global
potential-theoretic description of the bifurcation phenomena occurring within the quadratic polynomial family.
It is easy to obtain $\mu$ from $\Rc$ by pushing to $\P^2$ the slice of $\Rc$ by a suitable hypersurface of $X$. The measure $m$ can be obtained by pushing $\mu$ to the parameter space. Note that the support of $m$ is exactly the boundary of the Mandelbrot set. So $b\Mc$ is just the projection of the support of $\Rc$ to the parameter space.

Although  $\Rc$ is abstractly defined, our main result below shows that it can be obtained through a purely dynamical process which
makes it appearing as the Green $(2,2)$-current of the dynamical system $\widehat F :X\to X$. See \cite{DSCime} in order to compare with other dynamical systems. The support of $\Rc$ is then a kind of Julia set for $\widehat F$.

\begin{theorem}\label{t:main}
Let $\widehat F:X\to X$ be as above. Then, for any smooth closed $(2,2)$-form $\Omega$ on $X$ there exists a constant
 $\lambda_{\Omega}$ such that
 $$\lim_{n\to\infty} \frac{1}{n2^n}(\widehat F^n)^*(\Omega)=\lambda_{\Omega} \Rc$$
 in the sense of currents.
 Moreover, the constant $\lambda_\Omega$ only depends on the class of $\Omega$ in the Hodge cohomology group $H^{2,2}(X,\C)$.
\end{theorem}

We will see later that $\lambda_\Omega\not=0$ if and only if the class of $\Omega$ is outside a hyperplane of $H^{2,2}(X,\C)$. It is not difficult to see that $\lambda_\Omega$ depends linearly on $\Omega$ and $\lambda_\Omega\geq 0$ when $\Omega$ is a positive $(2,2)$-form.

We believe that the main ideas of our approach can be extended to families of holomorphic endomorphisms of the complex projective space in higher dimensions. This requires however  to solve several technical problems which are already non-trivial in the present setting.
The proof of the above result uses several tools from pluripotential theory and the theory of bifurcations. They will be presented in details in the next sections. 

In light of the above re-construction of the Mandelbrot set, we expect that the dynamical system $\widehat F :X\to X$ will give us more information about the bifurcations of the quadratic polynomial family. It is then useful to study the dynamics of $\widehat F$ systematically. More explicitly,  
the following open problem could be a first step, see \cite{Ah, DSCime, FJ, FS, Si, Ta} for similar equidistribution problems.

\begin{problem} \label{p:eq}
Let $\widehat F:X\to X$ be as above. Let $T$ be a  closed positive $(2,2)$-current on $X$. Find a sufficient condition on $T$ such that 
$$\lim_{n\to\infty} \frac{1}{n2^n}(\widehat F^n)^*(T)=\lambda_T \Rc$$
for some constant $\lambda_T>0$. 
\end{problem}

We expect that currents of integration on generic curves of $X$ satisfy the last equidistribution property. To end the introduction, we describe a heuristic argument which will be seen clearly in the proof of the main theorem and which is a starting point of our approach. It relates the Mandelbrot set with a phenomenon of curvature concentration. The system $F:\C^2\to\C^2$ permits to consider the dynamics of $p_c$
 in a family and the system $\widehat F:X\to X$ allows us to see how
 the system of $p_c$ varies as it takes into account the action of the
 differential of $F$.

For simplicity, consider a generic curve $Z$ of $\P^2$. Denote by $Y:=\widehat Z$ the lift of $Z$ to $X$ and 
$T:=[Y]$ the current of integration on $Y$. We have $(\widehat F^n)^*(T)=[\widehat F^{-n}(Y)]$ and the push-forward of $(\widehat F^n)^*(T)$ to $\P^2$ is equal to $[F^{-n}(Z)]$. It is not difficult to show that both $p_c^n$ and $F^n$ have degree $2^n$ while $\widehat F^n$ has degree $O(n2^n)$. So with the factor of normalization as in Theorem \ref{t:main} and
Problem \ref{p:eq}, we have 
$$\lim_{n\to\infty} \frac{1}{n2^n}[F^{-n}(Z)]= 0.$$
Therefore, the limit in Problem \ref{p:eq} is expected to be a vertical current which somehow represents the distribution of the curvature of $F^{-n}(Z)$ when $n$ goes to infinity. Over $\C\setminus b\Mc$, the dynamics of $p_c$ is structurally stable and the union of the most chaotic part $J_c$ with $c\not\in b\Mc$ is a union of a compact family of holomorphic graphs. We don't expect that the action of $F$ on $Z$ concentrate the curvature on the domain over 
$\C\setminus b\Mc$. Thus, the limit in Problem \ref{p:eq} will allow us to re-construct the Mandelbrot set.

\medskip

\noindent
{\bf Notation.} Throughout the paper, $p_c(z)=z^2+c$, $F(c,z)=(c,p_c(z))$, $X=\P\TP^2$, $\Pi:X\to \P^2$ and $\widehat F:X\to X$ are introduced above. In particular, $X$ is a $\P^1$-fibration over $\P^2$. Each point $q$ in $X$ corresponds to a complex tangent direction of $\P^2$ at the point $\Pi(q)$. For $(c,z)\in\C^2\subset \P^2$ and $v=(v_1,v_2)$ a non-zero tangent vector of $\C^2$ at $(c,z)$, the corresponding point of $X$ is denoted by $(c,z,[v])$ where $[v]=[v_1:v_2]$ is the homogeneous coordinate of the projective line $\Pi^{-1}(c,z)$.  If $Z$ is a complex curve in $\P^2$, we can lift it to a curve $Y:=\widehat Z$ in $X$ by taking the set of points $q$ such that $\Pi(q)$ is in $Z$ and the tangent direction of $\P^2$ at $\Pi(q)$, given by $q$, is also tangent to $Z$. The projection $\Pi:Y\to Z$ is then a finite map which is one-to-one outside the singularities of $Z$.

Denote by $L_\infty:=\P^2\setminus\C^2$ the projective line at infinity. Let 
$a_\infty$ (resp. $b_\infty$) be the point in $L_\infty$ of coordinates $[c:z]=[0:1]$ (resp. $[c:z]=[1:0]$).
Consider the hypersurface $V$ of the points $q$ in $X$ such that the tangent direction of $\P^2$, defined by $q$, is vertical. 
More precisely, $V$ is the union of the lifts of the projective lines through $a_\infty$.
Since $F$ preserves the vertical fibration of $\P^2$, the hypersurface $V$ is invariant under the action of $\widehat F$. For the affine coordinates $(c,z,t):=(c,z,v_1/v_2)$ of $X$, the hyperplane $V$ is of equation $t=0$. The projection $\Pi:V\to \P^2$ is just the blow-up of $\P^2$ at the point $a_\infty$ above.

Finally, we will use the standard Fubini-Study form $\omega_\FS$ on $\P^2$ normalized so that $\int_{\P^2} \omega_\FS^2=1$, and we fix a K\"ahler form $\omega_X$ on $X$ such that its restriction to each fiber of $\Pi$ has integral 1. They allow us to define the mass of a current and the volume of an analytic set in $\P^2$ or in $X$. Denote also by $\D$ and $\D(a,r)$ the unit Euclidean disc and the Euclidean disc of center $a$ and radius $r$ in $\C$. Denote for simplicity $\D_r:=\D(0,r)$. For any open subset $D$ of $\C$ and any holomorphic map $\gamma : D \to \C$ we denote by $\Gamma_\gamma$ the graph $\{(c,\gamma(c))\;\colon\;c\in D\}$ of $\gamma$ in $\C^2$. The notation $\ind_E$ stands for the characteristic function of a set $E$. The value of a current $T$ at a test form $\phi$ is denoted by $\langle T,\phi\rangle$ or $T(\phi)$. The notation $\{\cdot\}$ stands for the cohomology class of a closed current. 

\medskip

\noindent  \textbf{Acknowledgements:} The first author would like to thank the National University of Singapore for its support and hospitality during the visit where this work began.
The paper was partially written during the visit of the second author to the University of Cologne. He would like to thank this university, Alexander von Humboldt foundation and George Marinescu for their support and hospitality. The second author was also supported by the NUS grants C-146-000-047-001 and R-146-000-248-114.


\section{Properties of Green functions} \label{s:Green}

Recall that the dynamical Green function $g_c$ of $p_c$ is defined in \eqref{e:Green}. This is a continuous non-negative subharmonic function on $\C$ which is harmonic outside $J_c=bK_c$ and vanishes exactly on $K_c$. Its value $g_c(0)$ at 0 defines the Green function of the Mandelbrot set $\Mc$ which is continuous, non-negative, subharmonic, vanishing exactly on $\Mc$ and harmonic outside $b\Mc$, see \cite{Ber2} for details. 

In this section, we will give some properties of these Green functions that we will need later.
The proof of the following lemma is based on elementary potential theory on the complex plane and standard facts about dynamical stability 
within the quadratic family. We shall
in particular use the fact  that any connected component $\Omega$ of the interior of $\Mc$ is a stability component which either entirely consists of
hyperbolic parameters  (i.e. parameters $c$ for which the polynomial $p_c$ is hyperbolic) or entirely consists of non-hyperbolic ones. In the first case, the component is called \emph{hyperbolic} and 
the critical point $0$ belongs to the basin of some periodic attracting orbit of $p_c$ for every $c$ in $\Omega$. In the second case, the component is called \emph{non-hyperbolic}, the Fatou set of $p_c$
is reduced to the basin of $\infty$ and the critical point $0$ belongs to $J_c$ for every $c\in \Omega$. Conjecturally, such non-hyperbolic components do not exist.
See, for instance, \cite{Ber2}.

\begin{lemma}\label{l:Mandel-Green} 
Let $a$ and $b$ be two polynomials on $\C$ with $b$ not identically zero. Then the sequence of subharmonic functions 
$$\varphi_n(c):=\frac{1}{2^n}  \ln \vert b(c) p_c^n(0)-a(c)\vert$$ 
converges to the Green function $g_c(0)$ of $\Mc$ in $L^1_{\loc}(\C)$ when $n$ tends to infinity.
\end{lemma}

\proof 
First, observe that when $a$ and $b$ have a common factor, we can divide them by this factor because this operation doesn't change the $L^1_\loc$ limit of $\varphi_n$. So we can assume that $a$ and $b$ have no common factor.
We will use the following classical \emph {compactness
principle} : a sequence of subharmonic functions which is locally uniformly bounded from above and does not converge to $-\infty$ admits
a subsequence which is converging to some subharmonic function in $L^1_{\loc}$.

We first show that $\varphi_n(c)$ converges to $g_c(0)$ on $L^1_\loc(\C\setminus\Mc)$.  According to the above compactness principle, we only need to show that $\varphi_n(c)$ converges to $g_c(0)$ for almost every $c\in \C\setminus\Mc$. Fix a point $c\in \C\setminus\Mc$ such that $b(c)\not=0$.  We have seen at the beginning of the section that $g_c(0)\not=0$. It follows from \eqref{e:Green} that $p_c^n(0)$ tends to infinity as $n$ tends to infinity. Thus, using again \eqref{e:Green},
$$\lim_{n\to\infty} \varphi_n(c) = \lim_{n\to\infty}\frac{1}{2^n} \log|p_c^n(0)| +  \lim_{n\to\infty}\frac{1}{2^n} \log \Big|b(c)-{a(c)\over p_c^n(0)}\Big|=g_c(0)+0=g_c(0).$$
So we have  that $\varphi_n(c)$ converges to $g_c(0)$ on $L^1_\loc(\C\setminus\Mc)$.

We now want to prove that $(\varphi_n)_n$ is converging to $g_c(0)$ in $L^1_{\loc}(\C)$. By the compactness principle, 
this amounts to show that $g_c(0)$ is the only limit value of this sequence.  
Consider an arbitrary subsequence 
$(\varphi_{n_k})_k$  converging to a subharmonic function $\varphi$ in 
 $L^1_{\loc}(\C)$. We have shown that $\varphi=g_c(0)$ on $\C\setminus\Mc$ and to conclude that $\varphi =g_c(0)$ on whole $\C$, it is sufficient to show that $\varphi=0$ on $\Mc$. 
By the maximum principle, it is enough to check that  $\varphi =0$ on $b\Mc$ and
$\varphi$ is harmonic on $\Mc\setminus b\Mc$.
 
 It is known that the sequence $p_c^n(0)$ is bounded when $c$ is in $\Mc$. It follows from the definition of $\varphi_n$ that $\varphi\leq 0$ on $\Mc$. 
 On the other hand, since $\varphi(c)=g_c(0)> 0$ on $\C\setminus\Mc$ and $\varphi$ is upper semi-continuous as it is subharmonic, we deduce that $\varphi=0$ on $b\Mc$. It remains to establish that $\varphi$ is harmonic on any connected component $\Omega$ of the interior of $\Mc$. 
 We proceed by contradiction and assume that $\varphi\vert_{\Omega}$ is not harmonic.
  
Since $\Omega$ is a stable component, we may replace $(n_k)$ by a subsequence and assume  that $(p_c^{n_k}(0))_k$ converges locally uniformly to some holomorphic function $\sigma$ on $\Omega$. Take a disc $D\Subset \Omega$ such that $\ddc \varphi$ has positive mass on $D$.  Then, by definition of $\varphi$, the sequence of positive measures  
$$\frac{1}{2^{n_k}} \ddc \ln \vert b(c)p_c^{n_k}(0) -a(c)\vert $$ 
converges weakly to $\ddc \varphi$ in $D$. Note that $ \ddc \ln \vert b(c)p_c^{n_k}(0) -a(c)\vert$ is the sum of the Dirac masses at the zeros of 
$b(c)p_c^{n_k}(0) -a(c)$, counted with multiplicities.
Therefore, the number of zeros (counted with multiplicity) of $b(c)p_c^{n_k}(0) -a(c)$
 in $D$ tends to infinity as $k$ tends to infinity. Thus, by Hurwitz theorem,  $b(c)\sigma(c) -a(c)$ vanishes identically on $\Omega$. It follows that 
 $\sigma(c)=\frac{a(c)}{b(c)}$ on $\Omega$.
  
Observe that the curve $\Gamma_\sigma$ is not periodic (see the notation in the Introduction). Indeed, otherwise, by analytic continuation, 
  we would have   $p_c^N(\frac{a(c)}{b(c)})=\frac{a(c)}{b(c)}$ on $\C$ for some $N\geq 1$ which is clearly impossible for a degree reason. Since $\Gamma_\sigma$ is not periodic and $(p_c^{n_k}(0))_k$ converges locally uniformly to $\sigma$, we deduce that the stable component $\Omega$ is non-hyperbolic since otherwise
  the critical orbit should accumulate a periodic attracting cycle.
  
Let us pick a point $c_0\in \Omega$. As $\Omega$ is stable, there exists a dynamical holomorphic motion of the Julia sets centered at $c_0$ on $\Omega$. This is a continuous  map  
  $h: \Omega\times J_{c_0} \to \C$ of the form $(c,z)\mapsto h_c(z)$ which is one-to-one in $z$ when $c$ is fixed, holomorphic in $c$ when $z$ is fixed and satisfies the invariance relation
  $h_c(p_{c_0}(z))=p_c(h_c(z))$ for $(c,z)\in  \Omega\times J_{c_0}$. Since $\Omega$ is non-hyperbolic, the critical point $0$ belongs to $J_{c_0}$ and  $h_c(0)=0$ for every $c\in \Omega$, see Lemma \ref{l:NonHyp} below. 
Then we have  $p_c^{n_k}(0)=p_c^{n_k}(h_{c}(0))=h_c(p_{c_0}^{n_k}(0))$. It follows from the definition of $\sigma$ that $\sigma(c)=h_c(z_0)$ where $z_0:=\sigma(c_0)$.

Recall that $b(c)\sigma(c)  - a(c)=0$ on $\Omega$ and we have seen that $b(c')p_{c'}^{n_k}(0)-a(c')=0$ for some fixed $n_k$ large enough and some $c'\in\Omega$.
Recall also that $a$ and $b$ are assumed to have no common zero. So we must have $p_{c'}^{n_k}(0)=\sigma(c')$ which implies 
  $h_{c'} (p_{c_0}^{n_k}(0))=h_{c'}(z_0)$. Then, by the injectivity of the holomorphic motion, we obtain $p_{c_0}^{n_k}(0)=z_0$. Thus,
 $$p_c^{n_k}(0)= p_c^{n_k}(h_c(0))=h_{c} (p_{c_0}^{n_k}(0))=h_c(z_0)=\sigma(c)=\frac{a(c)}{b(c)}$$ 
 for all $c\in \Omega$. By analytic continuation, we obtain  $p_c^{n_k}(0)= \frac{a(c)}{b(c)}$ on $\C$ which is impossible as the left hand side is a polynomial of large degree in $c$. 
\endproof
  
  The following fact is well-known.
  
 \begin{lemma}\label{l:NonHyp}
 Let $c_0\in \Omega$ and $h: \Omega \times J_{c_0} \to \C$ be a holomorphic motion of Julia sets defined on some  non-hyperbolic component $\Omega$ of the interior of $\Mc$. Then $h_c(0)=0$ for every $c\in \Omega$.
 \end{lemma}
 
 \proof
As $\Omega$ is non-hyperbolic we know that  $0\in J_c$  for every $c\in \Omega$, moreover $h_{c_0}(0)=0$. We argue by contradiction and
 assume that $h_{c}(0)$ is not identically vanishing on $\Omega$.
 Let $(z_n)_n$ be a sequence of repelling periodic points of $p_{c_0}$ converging to $0$.
 Note that  $h_c(z_n)$ is a periodic repelling point of $p_c$ for every $c\in \Omega$.
Then, by  continuity of holomorphic motions, the sequence $(h_c(z_n))_n$ is locally uniformly
 converging 
 to $h_c(0)$ on $\Omega$ and, by Hurwitz theorem, there exist some $n$ big enough  and some $c_1\in \Omega$
such that $h_{c_1}(z_n)=0$. Thus, for $p_{c_1}$, the critical point $0$ belongs to some repelling
periodic cycle. This is clearly impossible.
\endproof
 
The following estimates are classical, see for instance \cite[Section 3.2.2]{BCime}. 
  
\begin{lemma}\label{l:Green}
The Green function $g_c(z)$ satisfies the following estimates on $\C\times\C$
\begin{enumerate}
\item[(1)] $g_c(z)\le\ln 2 + \max\left(\frac{1}{2} \ln \vert c\vert, \ln \vert z\vert\right)\; \text{when}\; \vert c\vert \ge 1\;;$
\item[(2)] $\max \left(g_c(z), \frac{1}{2} g_c(c)\right) \ge \ln \frac{\vert z\vert}{4}\; ;$
\item[(3)] for $10^{-2} \le \vert \alpha\vert  \le 10^{-1}$ and $\vert \beta \vert \le 2$, there exists $C_0\ge 1$ such that $g_c(\alpha c +\beta) < g_c(c)\;\textrm{when}\;   |c| \ge C_0$.
\end{enumerate}
\end{lemma}

\proof
(1) Choose $a$ such that $a^2=c$. Then, we have $|a|\geq 1$ and  for  $\vert z\vert \ge \vert a\vert$
$$\vert p_c(z)\vert= \vert z\vert^2  \Big\vert 1+\frac{c}{z^2}\Big\vert  \le 2\vert z\vert^2=2\max(\vert z\vert,\vert a \vert)^2.$$ 
By the maximum modulus principle, it follows that
$$\vert p_c(z)\vert\le 2\max(\vert z\vert,\vert a \vert)^2 \quad \text{for all} \quad (c,z)\in \C\times \C.$$
Since $\vert a\vert \ge 1$, by induction, one gets 
$$\vert p_c^n(z)\vert \le 2^{1+2+ \dots +2^{n-1}}  \max(\vert z\vert,\vert a\vert)^{2^{n}}$$
and this  implies Property (1).

\smallskip

(2) Consider the B\"{o}ttcher function of $p_c$
$$\varphi_c:\{g_c(z)>g_c(0)\}\to \P^1\setminus \overline{\D(0,e^{g_c(0)})}.$$ 
This is a univalent map satisfying $\varphi_c\circ p_c=\varphi_c^2$ and $\ln \vert \varphi_c\vert =g_c$. 
Denote by $\psi_c$ the inverse map of $\varphi_c$.
By Koebe $\frac{1}{4}$-theorem,  $\P^1\setminus \overline{\D(0,4r)}$ is contained in $\psi_c\big(\P^1\setminus \overline{\D(0,r)}\big)$
for every $r\geq e^{g_c(0)}$.

Pick $z\in \C$ and define $r:=\max(e^{g_c(z)},e^{g_c(0)})$. Then $z$ is outside $\psi_c\big(\P^1\setminus \overline{\D(0,r)}\big)$
since otherwise  $z=\psi_c(w)$ for some $w$ with $\vert w\vert >r$ and $e^{g_c(z)}=\vert \varphi_c(z)\vert=\vert w\vert >r\ge  e^{g_c(z)}$.
We thus have $\vert z\vert \le 4r=4\max(e^{g_c(z)},e^{g_c(0)})$ which, after taking logarithms, yields (2). We used here the invariance of $g_c$ which implies that 
$g_c(c)=g_c(p_c(0))=2g_c(0)$. 

\smallskip

(3) Take $C_0\ge 1$ sufficiently big so that $10^{-1} \vert c\vert +2 < 8^{-1} \vert c\vert$ and 
$10^{-2}\vert c\vert -2 \ge  \vert c\vert^{\frac{1}{2}}$ for $\vert c\vert \ge C_0$. Assume now that $\vert c\vert \ge C_0$.
Then $\vert \alpha c+ \beta \vert \ge 10^{-2} \vert c\vert -2  \ge \vert c\vert^{\frac{1}{2}}$ and, by (1), $g_c(\alpha c+\beta) \le \ln 2+ \ln \vert \alpha c+\beta\vert \le \ln 2 +  \ln (10^{-1} \vert c\vert +2) <\ln \frac{\vert c\vert}{4}$. The assertion (3) follows since, according to (2), we have
$\ln \frac{\vert c\vert}{4} \le g_c(c)$.
\endproof



\section{Some potential theoretic tools}\label{s:PTT}

We will give here some potential-theoretic results that we will use later.
These results are related to the following pre-order relation on probability measures. We refer to \cite{DSCime,Si, Tsu} for basic notions and properties of (pluri)potential theory. 

\begin{definition}\label{d:order}
Let  $\mu_1$ and $\mu_2$ be two positive measures  (e.g. two probability measures) 
with compact supports on a Stein manifold $M$. We set 
$\mu_1 \triangleright \mu_2$ and $\mu_2 \triangleleft \mu_1$ if $\langle \mu_1, \varphi\rangle \ge \langle \mu_2, \varphi\rangle$ for any plurisubharmonic (p.s.h. for short)
function $\varphi$ on $M$.
\end{definition}

We have the following property.

\begin{lemma} \label{l:order-mass}
Let $M$ and the measures $\mu_1,\mu_2$ be as in Definition \ref{d:order} with $\mu_1 \triangleright \mu_2$. Then we have  $\langle \mu_1, \varphi\rangle = \langle \mu_2, \varphi\rangle$ for any pluriharmonic function $\varphi$ on $M$. In particular, $\mu_1$ and $\mu_2$ have the same mass.
\end{lemma}
\proof
Since $\varphi$ is pluriharmonic, both $\varphi$ and $-\varphi$ are plurisubharmonic. We easily deduce the identity in the lemma from the definition of $\mu_1 \triangleright \mu_2$. Applying this identity to the function $\varphi=1$ implies that  $\mu_1$ and $\mu_2$ have the same mass.
\endproof

Recall that the equilibrium measure of a non-polar compact subset in $\C$ is characterized, among
 probability measures supported on this compact, as maximizing the energy. Using the relation $\triangleright$ this property can be rephrased as follows.

\begin{lemma}\label{l:Equi}
Let $K$ be a non-polar compact subset of $\C$ and $\mu_K$ be its equilibrium measure. Let $\nu$ be a probability measure supported on $K$.
If $\nu \triangleright \mu_K$ then $\nu=\mu_K$.
\end{lemma}

\proof 
Recall that the logarithmic potential and the energy of $\nu$ are defined by
$$v_\nu(z) :=\int \ln \vert z-a\vert \;d\nu (a) \quad \text{and} \quad I(\nu):=\iint \ln \vert z-a\vert d\nu(z)\;d\nu(a).$$ 
The ones for $\mu_K$ are defined in the same way. 
It is enough to show that $I(\nu)\ge I(\mu_K)$.

Recall that both $v_\nu$ and $v_{\mu_K}$ are subharmonic on $\C$. Hence,
it follows from Definition
\ref{d:order} and Fubini theorem that
\begin{eqnarray*}
I(\nu)=\int v_\nu(z)\;d\nu(z) \ge \int v_\nu(z)\;d\mu_K(z) = \int d\mu_K(z) \int \ln \vert z-a\vert\;d\nu(a) \qquad \\
=\int d\nu(a)  \int \ln \vert z-a\vert\;d\mu_K(z) \ge \int d\mu_K(a)  \int \ln \vert z-a\vert\;d\mu_K(z) =I(\mu_K). 
\end{eqnarray*}
This completes the proof of the lemma.
\endproof

Our aim is to extend the above lemma to some probability measures on $\C\times \C$. 
To this end we shall use the following unpublished result due to the second author and Sibony.

\begin{proposition}\label{p:DSU}
Let $K$ be a compact subset of $\C^n$. Let $\nu$ be a distribution supported by $K$. Then the following properties are equivalent:
\begin{enumerate}
\item[(i)] there exists a positive current $T$ of bi-dimension $(1,1)$ supported by $K$ such that $\ddc T=\nu$;
\item[(ii)] for every real-valued smooth function $\phi$ such that $\ddc\phi\geq 0$ on $K$, we have $\langle \nu,\phi\rangle \geq 0$.
\end{enumerate}
\end{proposition}

\proof
Assume  that $\ddc T=\nu$  for some positive current $T$ supported by $K$. Then,  for every smooth function $\phi$ satisfying $\ddc\phi\geq 0$ on $K$ we have 
$$\langle \nu,\phi\rangle=\langle \ddc T,\phi\rangle =\langle T, \ddc\phi\rangle \geq 0.$$
This shows that (i) $\Rightarrow$ (ii).

Let us now prove that 
(ii) $\Rightarrow$ (i).
 In the space of distributions, consider the convex cone
$$\Cc:=\big\{ \ddc T \;\colon\; \text{$T$ a positive current of bi-dimension $(1,1)$ supported by $K$} \big\}.$$
We first show that $\Cc$ is closed.

Assume that $T_n$ is a sequence of positive currents of bi-dimension $(1,1)$ supported by $K$ such that $\nu_n:=\ddc T_n$ converges to some distribution $\nu_\infty$. Using the K\"ahler form $\ddc\|z\|^2$ on $\C^n$, we have 
$$\|T_n\| =\langle T_n,\ddc\|z\|^2\rangle =\langle \nu_n, \|z\|^2\rangle \to \langle \nu_\infty,\|z\|^2\rangle.$$
So the mass of $T_n$ is bounded independently of $n$. Extracting a subsequence, we can assume that $T_n$ converges to a current $T$. Clearly, $T$ is positive and supported by $K$. So $\nu_\infty=\ddc T$ belongs to the cone $\Cc$.

Let us now show that if (i) is not true, that is $\nu \notin \Cc$,  then (ii) is not true either.  By Hahn-Banach theorem, there is a real-valued smooth function $\phi$ such that 
$$\langle \nu,\phi\rangle <  \langle \ddc T, \phi\rangle =\langle T, \ddc\phi\rangle $$
for all positive current $T$ of bi-dimension $(1,1)$ and supported by $K$. In particular, the inequality still holds if we multiply $T$ by any positive constant. When this constant tends to infinity, we see that $ \langle T, \ddc \phi\rangle\geq 0$. 
Since this is true for all $T$ positive supported by $K$, we deduce that $\ddc\phi \geq 0$ on $K$. 
For $T=0$, we get $\langle \nu,\phi\rangle<0$ and thus (ii) fails.
\endproof

We will state now the main result of this section. Let $m_1$ and $m_2$ be two probability measures with compact support in $\C$. 
For $m_i$-almost every $c\in\C$, $i=1,2$, consider a probability measure $\mu_{i,c}$ on a fixed compact subset $K$ of $\C\times\C$ whose support is contained in the vertical line $\{c\}\times\C$. Finally, define two probability measures on $K$ by
$$\mu_1:=\int  \mu_{1,c}\;dm_1(c) \quad \text{and} \quad \mu_2:=\int  \mu_{2,c}\;dm_2(c)$$
(we skip here the details about the dependence of $\mu_{i,c}$ on $c$ which is always assumed to be measurable).

\begin{lemma}\label{l:order-C2}
With the above notation, we assume moreover that 
$\mu_1 \triangleright \mu_2$. Then the following properties hold.
\begin{itemize}
\item[(1)]  We have $m_1 \triangleright m_2$;
\item[(2)] if moreover  $m_1 = m_2$, then $\mu_{1,c} \triangleright \mu_{2,c}$ for $m_1$-almost every $c\in \C$;
\item[(3)] if $\supp (m_1) \subset M$, $\supp (\mu_{1,c}) \subset K_c$ for some non-polar compact subsets $M,K_c$ of $\C$, and $m_2$, $\mu_{2,c}$ are the equilibrium measures of $M$ and $K_c$ respectively, then $\mu_1=\mu_2$.
\end{itemize}
\end{lemma}
\proof
The first assertion is immediately obtained from $\langle \mu_1,\varphi\rangle \ge  \langle \mu_2,\varphi\rangle$ for p.s.h. functions $\varphi$ on $\C^2$ which are only depending  on $c$. The third assertion follows from Lemma \ref{l:Equi} and the two former assertions.
It remains to establish the second assertion.

Without loss of generality, we can assume that the above compact set $K$ is a ball. Observe that any smooth function $\phi$ on $\C^2$ such that $\ddc\phi\geq 0$ on $K$ can be uniformly approximated by smooth p.s.h. functions on $\C^2$. Therefore, we can apply
Proposition \ref{p:DSU} to $\nu:=\mu_1-\mu_2$. So there exists a positive current $T$ of bi-dimension $(1,1)$ supported by $K$ 
such that
\begin{eqnarray}\label{e:Tmu} 
\ddc T=\mu_1-\mu_2.
\end{eqnarray}

Let  $\pi:\C^2\to\C$ be the canonical projection defined by $(c,z)\mapsto c$. Then $\pi_*(T)$ is a positive current with compact support.
Moreover, we have
$$\ddc \pi_*(T)=\pi_*(\ddc T)= \pi_*(\mu_1-\mu_2)=m_1-m_2=0.$$
So $\pi_*(T)$ is given by a constant function on $\C$ which, as  it has compact support, should be 0. We only need to consider the case $T\not=0$.

\medskip\noindent
{\bf Claim.}  There is a positive measure $\widetilde m$ with support in $\pi(K)$ and  positive currents $T_c$ with supports on 
$K\cap (\{c\}\times\C)$ for $\widetilde m$-almost every $c$, such that 
\begin{eqnarray}\label{e:Tm}
T=\int  T_c d\widetilde m(c).
\end{eqnarray}

\medskip

To prove the claim, we argue as in \cite[Lemma 3.3]{DSdens}. Let $\Cc$ be the convex cone of all  positive currents $S$ of bi-dimension $(1,1)$ supported in the compact set $K$ and satisfying $\pi_*(S) =0$. Let $\Cc_1$ be the subset of  $\Cc$ consisting of currents of  mass $1$. This is a compact convex set. Observing that $(\chi\circ \pi)S \in \Cc$ for any smooth positive function $\chi$ on $\C$ and every $S\in \Cc$,  
one sees that each extremal element of $\Cc_1$ is necessarily supported by $K\cap (\{c\}\times\C)$ for some $c$.
The decomposition \eqref{e:Tm} then follows from the classical Choquet's theorem.

\smallskip

Write for simplicity $m:=m_1=m_2$. From the assumptions and the identities \eqref{e:Tmu}, \eqref{e:Tm} we get 
\begin{eqnarray}\label{e:mu12}
\qquad \mu_1-\mu_2= \int  (\mu_{1,c}-\mu_{2,c}) dm(c) \quad \text{and} \quad \mu_1-\mu_2=\ddc T= \int  \ddc T_c d\widetilde m(c).
\end{eqnarray}
Observe that $\ddc T_c\not =0$  if and only if $T_c\not=0$ as $T_c$ has compact support in $\C^2$.

Write $\widetilde m=\widetilde m'+\widetilde m''$ where $\widetilde m'$ and  $\widetilde m''$ are respectively absolutely continuous and   singular with respect to $m$. 
Let $h\in L^1(m)$ be a non-negative function such that  $\widetilde m' = h m$.
After dividing  $\widetilde m'$ by $g:=h\ind_{h>0}$  and multiplying $T_c$ by $g(c)$ for $\widetilde m'$-almost every $c$ in the decomposition \eqref{e:Tm}, we may assume that $h$ only takes values $0$ or $1$ and we can set  $h=:\ind_E$ for some Borel subset $E$ of $\C$. It follows that $\widetilde m'=\ind_E\;m \leq m$ and the measure $m'':=m-\widetilde m'=\ind_{\C \setminus E}\;m$ is singular with respect to $\widetilde m'$. Define now $\widehat m:= \widetilde m'+m''+\widetilde m''=m+\widetilde m''=\widetilde m +m''$.

Since the measure $\widetilde m''$ is singular with respect to $m$, we may  define $\mu_{1,c}=\mu_{2,c}:=0$ for $\widetilde m''$-almost every $c$ and rewrite the first decomposition in \eqref{e:mu12} as
\begin{eqnarray*}
\mu_1-\mu_2= \int  (\mu_{1,c}-\mu_{2,c}) d\widehat m(c).
\end{eqnarray*}
Similarly, since $\widetilde m = \ind_E\;m + \widetilde m''$ 
 and $m''=\ind_{\C\setminus E}\;m$, we may set 
$T_c:=0$ for $m''$-almost every $c$ and rewrite the second decomposition in \eqref{e:mu12}  as
\begin{eqnarray*}
\mu_1-\mu_2= \int  \ddc T_c d\widehat m(c).
\end{eqnarray*}
Defining $R_c:=(\mu_{1,c}-\mu_{2,c})-\ddc T_c$ for $\widehat m$-almost every $c$ we thus get
$$\int  R_c d\widehat m(c)=0.$$ 

Consider a smooth test function $\phi$ on $\C^2$ and define $l(c):=\langle R_c,\phi\rangle$. If $\chi$ is any smooth test function in $c$, the above identity applied to the test function $\chi(c)\phi(c,z)$ gives
$$\int  \chi(c) l(c)d\widehat m(c)=0.$$
Since this is true for every $\chi$, we get  $l(c)=0$ for $\widehat m$-almost every $c$. 

Using a countable dense family of test functions $\phi$, we deduce that 
$R_c=0$ for $\widehat m$-almost every $c$ which, by
Proposition \ref{p:DSU}, yields $\langle \mu_{1,c} - \mu_{2,c},\phi\rangle\geq 0$ for all smooth p.s.h. functions on $\C^2$. The same property holds for any p.s.h. function on $\C^2$ because we can approximate it by a decreasing sequence of smooth p.s.h. ones.
Thus, $\mu_{1,c}\triangleright \mu_{2,c}$ and this completes the proof of the lemma. 
\endproof

We present now a situation where the above pre-order of probability measures naturally appears and is useful. We refer to \cite{DSU} for some details.

Let $M\subset \C^m$ and $N\subset \C^n$ be two bounded open sets. Denote by $\pi_M:M\times N\to M$ and $\pi_N:M\times N\to N$ the canonical projections. Consider two  closed positive currents $T$ and $S$ on $M\times N$ of bi-dimension $(m,m)$ and $(n,n)$ respectively.
Assume that $T$ is {\it horizontal-like} in the sense that $\pi_N(\supp(T))$ is relatively compact in $N$. Similarly, we will say that  $S$ is {\it vertical-like} if $\pi_M(\supp (S))$ is relatively compact in $M$.

The current $(\pi_M)_*(T)$ is a  closed positive current of bi-dimension $(m,m)$ of $M$. So it is defined by a non-negative constant function that we denote by $\lambda_T$ and call  {\it the vertical mass} of $T$. The following result is a consequence of \cite[Prop. 3.3]{DSU}.

\begin{lemma} \label{l:mass-horiz}
Let $K$ be a compact subset of $M$. Then there is a constant $c>0$ independent of $T$ and $\lambda_T$ such that the mass of $T$ on $K\times\C^n$ is bounded by $c\lambda_T$. 
\end{lemma}

The slice of $T$ by the space $\{a\}\times\C^n$ is well defined for every $a\in M$. We denote it by $\langle T,\pi_M,a\rangle$. This is a positive measure of mass $\lambda_T$ with compact support in $\{a\}\times N$. It can be obtained in the following way.

Fix any smooth non-negative radial function $\psi$ with compact support in $\C^m$ with integral 1. Define for $\epsilon>0$ the function 
$\psi_{\epsilon,a}(z) := \epsilon^{-2m} \psi (\epsilon^{-1}(z-a))$ which approximates the Dirac mass at $a$ when $\epsilon$ goes to 0. 
Then, for any given smooth function $\phi$ on $\C^m\times \C^n$, 
one has
$$\langle T,\pi_M, a\rangle (\phi)=\lim_{\epsilon\to 0} \big\langle T\wedge (\pi_M)^* (\psi_{\epsilon,a} \Theta_m),\phi \big\rangle,$$
where $\Theta_m$ is the standard volume form of $\C^m$.
Moreover, when $\phi$ is p.s.h., the functions 
\begin{equation} \label{e:Phi}
\Phi_\epsilon(a):=\big\langle T\wedge (\pi_M)^* (\psi_{\epsilon,a} \Theta_m),\phi \big\rangle \quad \text{and} \quad  
 \Phi(a):= \langle T,\pi_M, a\rangle (\phi)
 \end{equation} 
 are p.s.h. on $M$ and $\Phi_\epsilon$ decreases to $\Phi$ as $\epsilon$ decreases to 0, see \cite{DSU}.

\begin{lemma}\label {l:slice-order}
Let $(T_n)_n$ be a sequence of closed positive horizontal-like currents of bi-dimension $(m,m)$ on $M \times N$ converging to some  closed positive horizontal-like current. Let $a$ be a point in $M$.
Assume that the sequence of measures $(\langle T_n,\pi_M, a\rangle)_n$ is also convergent. Then we have the following property on $\C^m\times\C^n$
$$\big\langle \lim_{n\to\infty} T_n,\pi_M,a \big\rangle \triangleright \lim_{n\to\infty} \langle  T_n,\pi_M,a\rangle.$$
\end{lemma}

\proof 
Fix a smooth p.s.h. function $\phi$ on $\C^m\times\C^n$. 
Let $T$ be the limit of $T_n$ which is a closed positive horizontal-like current. Define $\Phi_\epsilon$ and $\Phi$ as in \eqref{e:Phi}. Denote 
the analogous functions associated to $T_n$ by $\Phi_{n,\epsilon}$ and $\Phi_n$.  

Fix an $\epsilon>0$. According to the above discussion on the slice of horizontal-like currents, we have 
$\Phi_n(a)\leq \Phi_{n,\epsilon}(a)$. Since $T_n$ converges to $T$, we deduce that $\Phi_{n,\epsilon}(a)$ tends to $\Phi_\epsilon(a)$ as $n$ tends to infinity. It follows that 
$\lim \Phi_n(a) \leq \Phi_\epsilon(a)$. Taking $\epsilon$ going to 0 gives $\lim \Phi_n(a) \leq \Phi(a)$.
Equivalently, we have
$$\big\langle \lim_{n\to\infty} T_n,\pi_M,a \big\rangle (\phi) \geq \lim_{n\to\infty} \langle  T_n,\pi_M,a\rangle (\phi).$$

This property still holds for any p.s.h. function $\phi$ because we can approximate it by a decreasing sequence of smooth ones.
The lemma follows.
\endproof

\begin{remark} \label{r:vertical} \rm
Clearly, we can apply Lemmas \ref{l:mass-horiz} and \ref{l:slice-order} to vertical-like currents. Indeed, the involution $(z,w)\mapsto (w,z)$ transforms vertical-like currents to horizontal-like ones.
\end{remark}


\section{Cohomological arguments}\label {s:coh}

Our aim in this section is to show, using the cohomology of $X$, that the proof of Theorem \ref{t:main} can be reduced to the case where $\Omega$ is a special form. 

\smallskip\noindent
{\bf Some basic properties of the map $\widehat F$}.
We use the affine coordinates $(c,z)$ of $\P^2$ and $(c,z,t)$ of $X$ given in the Introduction with $t:=v_1/v_2$. It is not difficult to see that 
\begin{equation} \label{e:F-hat}
F^n(c,z)=(c,p_c^n(z)) \quad \text{and}\quad \widehat F^n(c,z,t)=\Big(c,p_c^n(z), {t+\partial p_c^n(z)/\partial c \over \partial p_c^n(z)/\partial z} \Big).
\end{equation}
So $\widehat F^n$ is a rational map and it is holomorphic in some Zariski open subset $U_n$ of $X$. 
Its topological degree, i.e. the number of points in a generic fiber of $\widehat F^n$, is equal to the one of $F^n$ and hence is equal to $2^n$. 

Denote by $\Gamma_n$ the closure of the graph of $\widehat F^n$ on $U_n$. This is an irreducible analytic subset of dimension 3 of 
$X\times X$ that we also call {\it the graph of $\widehat F^n$ on $X$}. It doesn't depend on the choice of $U_n$. 
Denote by $\pi_j:X\times X\to X$, with $j=1,2$,  the canonical projections. 
The two {\it indeterminacy sets} of $\widehat F^n$ are defined by
$$I_{j,n}:=\big\{q\in X,\ \dim \pi_j^{-1}(q)\cap\Gamma_n \geq 1\big\} \quad \text{for} \quad j=1,2.$$
These are analytic subsets of dimension at most 1 of $X$ which play an important role in the study of the action of $\widehat F$ on currents.
Observe that $\pi_1:\Gamma_n\to X$ is 1:1 outside the analytic set $\pi_1^{-1}(I_{1,n})$. 

In order to better understand the action of $\widehat F^n$ on cohomology, we need the following property of $\Gamma_n$. We refer to the end of the Introduction for the definition of $V$ and $L_\infty$.

\begin{lemma} \label{l:I2}
The intersection $I_{2,n}\cap V$ is contained in $\Pi^{-1}(L_\infty)$.
\end{lemma}
\proof
Consider  a point $q\in V$ such that $\Pi(q)=(c,z)\in\C^2$. We have to show that $q$ is not in $I_{2,n}$. 
For this purpose, it is enough to show that the set  $A:=\pi_2^{-1}(q)\cap\Gamma_n$ is finite. Observe that $\pi_1$ is injective on $A$ because $A$ is contained in $X\times \{q\}$. So we only need to check that $A_0:=\pi_1(A)=\widehat F^{-n}(q)$ is finite.

Consider a point $q_0\in A_0$ and define $(c_0,z_0):=\Pi(q_0)$. We necessarily have $F^n(c_0,z_0)=(c,z)$ or equivalently $c_0=c$ and $z_0\in p_c^{-n}(z)$. 
There are finitely many $(c_0,z_0)$ satisfying these properties. 
It remains to show that there are finitely many  tangent directions $[v_0]$ at $(c_0,z_0)$ which are sent by the differential $dF^n(c_0,z_0)$ to the vertical direction or to 0. Since $F^n$ preserves the vertical lines of $\C^2$, it is clear that only the vertical direction can satisfy the last property. The lemma follows.
\endproof

\smallskip\noindent
{\bf Action on smooth forms and cohomology.} 
Let $[\Gamma_n]$ denote the current of integration on $\Gamma_n$. The pull-back action of $\widehat F^n$ on a current $T$ is defined by
$$(\widehat F^n)^*(T):= (\pi_1)_*(\pi_2^*(T)\wedge [\Gamma_n])$$
when the last wedge-product is well defined. We consider here a particular case. 

Assume that $T$ is given by a smooth differential form.
Then $(\widehat F^n)^*(T)$ is well-defined and given by a differential form with $L^1$ coefficients. So it has no mass on proper analytic subsets of $X$. If moreover $T$ is closed, exact or positive, so is  $(\widehat F^n)^*(T)$. Therefore, the above pull-back operator defines linear actions on the Hodge cohomology groups that we still denote by the same notation
$$(\widehat F^n)^* : H^{p,p}(X,\C) \to H^{p,p}(X,\C) \quad \text{for} \quad p=0,1,2,3.$$
As for a general rational map, the last operator is identity when $p=0$ and is the multiplication by the topological degree when $p=3$.
We are interested now in the case where $p$ is 1 or 2.

Recall that $\omega_X$ is a K\"ahler form on $X$ fixed at the end of the Introduction. So $\omega_X^p$ is a positive $(p,p)$-form and $(\widehat F^n)^*(\omega_X^p)$ is a closed positive $(p,p)$-current. As for all closed positive currents, its mass only depends on its cohomology class in $H^{p,p}(X,\C)$ and is given by
$$\|(\widehat F^n)^*(\omega_X^p)\|:=\big\langle (\widehat F^n)^*(\omega_X^p), \omega_X^{3-p}\big\rangle.$$
It was proved in \cite{DSEntropy} that the norm of $(\widehat F^n)^*$ on $H^{p,p}(X,\C)$ satisfies 
\begin{equation} \label{e:norm-mass}
A^{-1} \|(\widehat F^n)^*(\omega_X^p)\| \leq \|(\widehat F^n)^*\|_{H^{p,p}} \leq A \|(\widehat F^n)^*(\omega_X^p)\|
\end{equation}
for some constant $A\geq 1$ independent of $\widehat F$ and of $n$. 

\begin{lemma} \label{l:coh-11}
We have $\|(\widehat F^n)^*\|_{H^{1,1}}= O(2^n)$ when $n$ tends to infinity.
\end{lemma}
\proof
It is known that the behavior of $\|(\widehat F^n)^*\|_{H^{1,1}}$ doesn't change when we use a birational modification of $X$. More precisely, if $\pi:X\to \widetilde X$ is a birational map from $X$ to another projective threefold $\widetilde X$ and $\widetilde F^n:=(\pi\circ \widehat F \circ \pi^{-1})^n$ is conjugated to $\widehat F^n$ by $\pi$, then 
$$A^{-1} \|(\widetilde F^n)^*\|_{H^{1,1}}  \leq \|(\widehat F^n)^*\|_{H^{1,1}} \leq A \|(\widetilde F^n)^*\|_{H^{1,1}} $$
for some constant $A\geq 1$ independent of $n$, see \cite{DSEntropy}.  

From \eqref{e:F-hat}, we have a rational map on $\C^3$ that extends to $\widehat F^n:X\to X$. We can also extend it to a rational map $\widetilde F^n:\P^3\to\P^3$. Since $\dim H^{1,1}(\P^3,\C)=1$, the action of $\widetilde F^n$ on $H^{1,1}(\P^3,\C)$ is just the multiplication by some positive constant $\lambda_n$. If $H$ is a generic hyperplane in $\P^3$, the formula in \eqref{e:F-hat} implies that $\widetilde F^{-n}(H)$ is a hypersurface of degree $O(2^n)$. It follows that $\lambda_n=O(2^n)$ which ends the proof of the lemma.
\endproof

The following two lemmas, together with the fact that the dimension of $H^{2,2}(X,\C)$ is $2$, will allow us to prove Theorem \ref{t:main} by only considering a suitable form $\Omega$ (see also Lemma \ref{l:coh-22} and its proof below).

\begin{lemma} \label{l:coh-only}
Let $\Omega$ and $\Omega'$ be two smooth real closed $(2,2)$-forms on $X$ having the same cohomology class in $H^{2,2}(X,\C)$. Then
$$\lim_{n\to\infty} \Big[{1\over n2^n} (\widehat F^n)^*(\Omega) -  {1\over n2^n} (\widehat F^n)^*(\Omega') \Big]=0$$
in the sense of currents.
\end{lemma}
\proof
By the classical $\ddc$-lemma, there is a smooth real $(1,1)$-form $\alpha$ on $X$ such that $\Omega-\Omega'=\ddc \alpha$. Adding to $\alpha$ a constant times $\omega_X$ we can assume that $\alpha$ is positive. We can also divide $\Omega, \Omega'$ and $\alpha$ by a constant and assume that $\alpha\leq \omega_X$.

The expression in the brackets in the lemma is equal to ${1\over n2^n} \ddc (\widehat F^n)^*(\alpha)$. So in order to obtain the lemma, it is enough to show that $\| (\widehat F^n)^*(\alpha)\|=O(2^n)$. Moreover, since $\alpha\leq \omega_X$, it is enough to check that $\| (\widehat F^n)^*(\omega_X)\|=O(2^n)$. By \eqref{e:norm-mass}, we only need to prove that $\|(\widehat F^n)^*\|_{H^{1,1}}=O(2^n)$. But this is given by Lemma \ref{l:coh-11} above. 
\endproof

We fix a smooth positive $(2,2)$-form $\alpha_0$ on $\P^2$ with integral 1. Define $\Omega_0:=\Pi^*(\alpha_0)$. This is a positive closed $(2,2)$-form on $X$. We have the following lemma.

\begin{lemma} \label{l:Omega-0}
The sequence ${1\over n2^n} (\widehat F^n)^*(\Omega_0)$ tends to $0$ in the sense of currents when $n$ goes to infinity.
\end{lemma}
\proof
By definition of $\widehat F$ and $\Omega_0$, we have 
$${1\over n2^n} (\widehat F^n)^*(\Omega_0) = \Pi^*\Big({1\over n2^n} (F^n)^*(\alpha_0) \Big).$$
Now, observe that $\alpha_0$ defines a probability measure on $\P^2$ and the action of $F^n$ on a positive measure multiplies its mass by the topological degree $2^n$ of $F^n$. We deduce that ${1\over n2^n} (F^n)^*(\alpha_0)$ is a positive measure of mass $1/n$. It is clear that this measure tends to 0 when $n$ goes to infinity. The lemma follows.
\endproof

\smallskip\noindent
{\bf Action on positive closed currents and cohomology.} Let $T$ be a positive closed $(p,p)$-current on $X$. Assume that $T$ vanishes in a neighbourhood of $I_{2,n}$. Then the current $(\widehat F^n)^*(T)$ is well-defined, see \cite{DSPull}. Moreover, the action of $\widehat F^n$ on $T$ is compatible with the action on cohomology, i.e. we have 
$$\{(\widehat F^n)^*(T)\} = (\widehat F^n)^*\{T\}.$$

In particular, for $p=2$, if $T$ is the current of integration on a generic analytic curve of $X$ then $(\widehat F^n)^*(T)$ is well-defined. Indeed, since $I_{2,n}$ has dimension at most 1, a generic analytic curve in $X$ has no intersection with $I_{2,n}$. 

\begin{lemma} \label{l:coh-22}
Let $\theta$ be a cohomology class in $H^{2,2}(X,\C)$. Denote by $\lambda$ the complex number such that the class of $\Pi_*(\theta)$ in $H^{1,1}(\P^2,\C)$ is $\lambda$ times the class of a projective line. Then the class ${1\over n2^n} (\widehat F^n)^*(\theta)$ converges to ${\lambda\over 2}\{\Omega_0\}$  when $n$ goes to infinity.
\end{lemma}
\proof
Recall that $H^{p,q}(\P^2,\C)=0$ for $p\not=q$ and $\dim H^{p,p}(\P^2,\C)=1$ for $p=0,1,2$. It follows from Leray's spectral theory that 
$\dim H^{p,p}(X,\C)=2$ for $p=1,2$, see e.g. \cite[Th.\,7.33]{Voisin}.
Fix any class $\theta_0$ in $H^{2,2}(X,\C)$ such that $\Pi_*(\theta_0)$ is equal to the class of a projective line in $H^{1,1}(\P^2,\C)$. 
Since $\Pi_*(\Omega_0)=0$, we deduce that $\{\Omega_0\}$ and $\theta_0$ constitute a basis of $H^{2,2}(X,\C)$. 
Therefore, the class $\theta-\lambda\theta_0$ is co-linear to $\{\Omega_0\}$. By Lemma \ref{l:Omega-0}, the class ${1\over n2^n} (\widehat F^n)^*(\theta-\lambda\theta_0)$ tends to 0 as $n$ tends to infinity. Thus, we only need to prove the lemma for $\theta_0$ instead of $\theta$. 

Consider a generic projective line $L$. We can choose $\theta_0$ as the cohomology class of the current $[\widehat L]$. 
We first prove the following claim.

\medskip\noindent
{\bf Claim 1.}  Let $\eta$ be any limit value of ${1\over n2^n} (\widehat F^n)^*(\theta_0)$. Then $\eta$ is co-linear to $\{\Omega_0\}$.

\medskip

Since $\{\Omega_0\}$ and $\theta_0$ constitute a basis of $H^{2,2}(X,\C)$, we only need to show that $\Pi_*(\eta)=0$. This is clear because using that $\Pi_*(\{\widehat L\})=\{L\}$, we have 
$$\lim_{n\to\infty} \Pi_*\Big({1\over n2^n} (\widehat F^n)^*(\theta_0)\Big) = \lim_{n\to\infty} {1\over n2^n} (F^n)^*\{L\}=0.$$
We used here the fact that the action of $(F^n)^*$ on $H^{1,1}(\P^2,\C)$ is the multiplication by $2^n$. So Claim 1 is true. 

Now, observe that $\{\Omega_0\}\smallsmile \{V\}=1$ because $\pi :V\cap\pi^{-1}(\C^2)\to\C^2$ is bi-holomorphic. So, in order to get the lemma, it is enough to show that the number of points in the intersection $\widehat F^{-n}(\widehat L)\cap V=\widehat{F^{-n}(L)}\cap V$, counted with multiplicities, is equal to $n2^{n-1}+O(2^n)$. 
This is a consequence of the following two claims.

\medskip\noindent
{\bf Claim 2.} The number of points of the intersection $\widehat F^{-n}(\widehat L)\cap V$ in $\pi^{-1}(\C^2)$, counted with multiplicities, is equal to $n2^{n-1}$ for $n\geq 1$.
\medskip

By Lemma \ref{l:contact-lift} below, this number is equal to the number of points in the intersection between the critical set of $F^n$ and $F^{-n}(L)$ in $\C^2$.
Denote by $C$ the line $z=0$ which is the critical set of $F$. Then the critical set of $F^n$ is the union of the curves $C, F^{-1}(C),\ldots, F^{-n+1}(C)$. We now count the number of points of $F^{-m}(C)\cap F^{-n}(L)$ in $\C^2$ for $0\leq m\leq n-1$. By taking the image by $F^m$, we see that this number is equal to the number of points in the intersection of $F^m(F^{-m}(C))$ and $F^{-n+m}(L)$, counted with multiplicities. Observe that $F^m(F^{-m}(C))$ is equal to $C$ with multiplicity $2^m$ and if $L$ has equation $z=\alpha c+\beta$ for $\alpha,\beta\in\C$, then the equation of $F^{-n+m}(L)$ is $p_c^{n-m}(z)=\alpha c+\beta$. Its intersection with $C$ is given by the solutions of the equation $p_c^{n-m}(0)=\alpha c+\beta$. Since the left hand side is a polynomial of degree $2^{n-m-1}$ in $c$, the last equation has $2^{n-m-1}$ solutions counting multiplicity. We conclude that the number of points of $F^{-m}(C)\cap F^{-n}(L)$ in $\C^2$ is $2^m\times 2^{n-m-1}=2^{n-1}$. This implies Claim 2.

\medskip\noindent
{\bf Claim 3.} The number of points of the intersection $\widehat F^{-n}(\widehat L)\cap V$ in $\pi^{-1}(L_\infty)$, counted with multiplicities, is at most equal to $2^n-1$ for $n\geq 1$.
\medskip

The equation of $F^{-n}(L)$ is $p_c^n(z)=\alpha c+\beta$. Since $p_c^n(z)$ is a polynomial of degree $2^n$ in $(z,c)$ whose unique highest degree term is $z^{2^n}$, the curve $F^{-n}(L)$ intersects $L_\infty$ at a unique point $b_\infty=[1:0]$, see the end of the Introduction for the notation. 
It follows from the definition of $V$ that if the intersection $\widehat F^{-n}(\widehat L)\cap V\cap\pi^{-1}(L_\infty)$ is non-empty, it should be the singleton $b_\infty':=(b_\infty,[v_\infty])$, where $v_\infty$ is a tangent vector of $L_\infty$ at $b_\infty$. We need to estimate the multiplicity of $\widehat F^{-n}(\widehat L)\cap V$ at this point $b_\infty'$.
 
We will use the local coordinates $c':=1/c$ and  $z':=z/c$ of $\P^2$ near $b_\infty$. So $c=1/c'$, $z=z'/c'$ and $b_\infty=(0,0)$ in these coordinates. Define for simplicity $q(c,z):=p_c^n(z)-\alpha c-\beta$ and $r(c',z'):=c'^{2^n}q(c,z)$ which is a polynomial in $c',z'$ whose zero set is $F^{-n}(L)$. So the intersection between $F^{-n}(L)$ and $L_\infty$ is given by the equations $r(c',z')=0$ and $c'=0$. If we replace $c'$ by $0$, the equation $r(c',z')=0$ becomes $z'^{2^n}=0$. Thus, $b_\infty$ is a point of intersection of order $2^n$ between $F^{-n}(L)$ and $L_\infty$.

Let $S$ be any irreducible germ of $F^{-n}(L)$ at $b_\infty$ and let $m$ denote the multiplicity of its intersection with $L_\infty$ at $b_\infty$. We will show that the intersection between the lift $\widehat S$ of $S$ to $X$ and the hypersurface $V$ at the point $b_\infty'$ is smaller than $m$. This implies Claim 3. We will use the local coordinates $c',z'$ and $t':=v_1'/v_2'$ of $X$ at $b_\infty'$ such that $v'=(v_1',v_2')$ is a tangent vector of $\P^2$ at the point $(c',z')$. In these coordinates, we have $b_\infty'=(0,0,0)$ and $V$ is given by $t'=0$.

We can parametrize the curve $S$ by $(s^p, s^qh(s))$ for $p,q\geq 1$, $s\in \C$ small and $h(s)$ a non-vanishing holomorphic function. The intersection between $S$ and $L_\infty$ is given by $s^p=0$. We deduce that $p=m$. If $S$ is not tangent to $L_\infty$ at $b_\infty$, then $\widehat S$ does not contain $b_\infty'$ and we have the desired property. Otherwise, we have $q<p=m$. It is not difficult to see that the curve 
$\widehat S$ is parametrized by $(s^p,s^qh(s), s^{p-q}l(s))$ for some holomorphic function $l(s)$ with $l(0)\not=0$. Its intersection with $V$ is given by the equation $s^{p-q}l(s)=0$ and hence $b_\infty'$ is an intersection point of multiplicity $p-q<m$. This ends the proof of Claim 3 and  the proof of the lemma as well.
\endproof

\begin{lemma} \label{l:contact}
Let $L$ be a generic projective line in $\P^2$. Then the curve $F^{-n}(L)$ has no singularity in $\C^2$ and the order of contact of any vertical line $\{c\}\times\C$ with $F^{-n}(L)$ at each point is equal to $0$ or $1$. Moreover, a vertical line $\{c\}\times\C$ is tangent to $F^{-n}(L)$ at a point $(c,z)$ if and only if $(c,z)$ belongs to the critical set of $F^n$.  
\end{lemma}
\proof
Observe that the critical set of $F^n$ in $\C^2$ is given by the equation $\partial p^n_c(z)/\partial z=0$ and its image by $F^n$ is the set of critical values of $F^n$ in $\C^2$. We only consider a generic projective line $L$ such that
\begin{itemize}
\item $L$ is not a vertical line;
\item  $L$ intersects transversally the set of critical values of $F^n$ in $\C^2$; in particular, $L$ contains no singular point of 
the set of critical values of $F^n$ in $\C^2$;
\item if $z\in\C$ is a multiple zero of $\partial p^n_c(z)/\partial z$ when we fix $c\in\C$, then $F^n(c,z)$ does not belong to $L$; in particular, near each point of its intersection with $F^{-n}(L)$ in $\C^2$, the critical set of $F^n$ is a holomorphic graph over an open set of the $c$-axis of $\C^2$. 
\end{itemize}

Fix a point $(c_0,z_0)$ in $F^{-n}(L)\cap\C^2$ such that  $\partial p_c^n(z)/\partial z$ does not vanish at $(c_0,z_0)$. Then $F^n$ defines a local bi-holomorphism between a neighbourhood of $(c_0,z_0)$ and a neighbourhood of $F^n(c_0,z_0)$. Hence $F^{-n}(L)$ is smooth at $(c_0,z_0)$. Moreover, 
since $F$ preserves the vertical fibration of $\C^2$, it is easy to see that $F^{-n}(L)$ is not tangent to the line $\{c_0\}\times \C$ at  $(c_0,z_0)$. 
So the vertical line $\{c_0\}\times\C$ intersects $F^{-n}(L)$ transversally at the point $(c_0,z_0)$.

Consider now a point  $(c_0,z_0)$ in $F^{-n}(L)\cap\C^2$ such that  $\partial p_c^n(z)/\partial z$ vanishes at $(c_0,z_0)$. So $(c_0,z_0)$ is a critical point of $F^n$. Since $L$ is generic as described above, $z_0$ is a simple zero of the polynomial $\partial p_{c_0}^n(z)/\partial z$ and the critical set of $F^n$ near $(c_0,z_0)$ is a graph of some holomorphic function $h(c)$ over a neighbourhood of $c_0$ in the $c$-axis of $\C^2$.
We now use the local coordinate system $c':=c-c_0$, $z':=z-h(c)$ near the point $(c_0,z_0)$ and the local coordinate system $c'':=c-c_0$, $z''=z-p_c^n(h(c))$ near the point $F^n(c_0,z_0)$. In these coordinates, the critical set and the set of critical values of $F^n$ are given by $z'=0$ and $z''=0$ respectively. We only work near the point $(0,0)$.

So we see that $F^n$ has the form
$(c',z')\mapsto (c', z'^2g(c',z'))$, where $g(c',z')$ is a non-vanishing holomorphic function. The line $L$ has the form $c''=z''l(z'')$ for some holomorphic function $l(z'')$ which does not vanish because $L$ is transverse to the critical values of $F^n$.
We see that $F^{-n}(L)$ is given by an equation $c'=z'^2\widetilde g(c',z')$ for some non-vanishing holomorphic function $\widetilde g$. It is now clear that the contact order between $F^{-n}(L)$ and the vertical 
line $\{c'=0\}$ at the point $(0,0)$ is equal to 1. We also see that $F^{-1}(L)$ is smooth at the point $(0,0)$, or equivalently, at the point $(c_0,z_0)$ in the original coordinates. This ends the proof of the lemma.
\endproof

We deduce from the last lemma and the definition of $V$ the following result.

\begin{lemma} \label{l:contact-lift}
Let $L$ be a generic projective line in $\P^2$. Then $\widehat F^{-n}(\widehat L)$ intersects $V\cap \pi^{-1}(\C^2)$ transversally. Moreover, this intersection is exactly the set of points $(c,z,[v])$ in $\pi^{-1}(\C^2)$ such that $(c,z)$ is an intersection point between $F^{-n}(L)$ and the critical set of $F^n$, and $[v]$ is the tangent direction of  $F^{-n}(L)$ at $(c,z)$ which is also the vertical direction.  
\end{lemma}

\begin{proposition} \label{p:limit-22}
Let $\Omega$ be a smooth closed positive $(2,2)$-form on $X$. Let $\lambda$ be the mass of the current $\Pi_*(\Omega)$ on $\P^2$. Then the mass of  ${1\over n2^n} (\widehat F^n)^*(\Omega)$ is bounded independently of $n$. Moreover, if $T$ is any limit value of   ${1\over n2^n} (\widehat F^n)^*(\Omega)$, then there is a positive measure $\nu$ of mass ${\lambda \over 2}$ on $\P^2$ such that $T=\Pi^*(\nu)$.
\end{proposition}
\proof
Recall that the mass of a positive closed current depends only on its cohomology class. Therefore, the first assertion is a direct consequence of Lemma \ref{l:coh-22}. This lemma also implies that when $T=\Pi^*(\nu)$ for some positive measure $\nu$ on $\P^2$, then the mass of $\nu$ is equal to ${\lambda \over 2}$. So, it remains to prove the existence of $\nu$ such that $T=\Pi^*(\nu)$. In other words, the current $T$ is vertical in the sense of  \cite{DSdens}. It was shown in this reference that $T$ is vertical if and only if $T\wedge \Pi^*(\omega_\FS)=0$. 
The last identity is clear because  $T\wedge \Pi^*(\omega_\FS)$ is a positive measure and its cohomology class is equal to 
${\lambda \over 2} \{\Omega_0\} \smallsmile \Pi^*\{\omega_\FS\}=0$, according to Lemma \ref{l:coh-22}. This ends the proof of the proposition. 
\endproof

\section{Proof of the main result} \label{s:proof}

We will give in this section the proof of Theorem \ref{t:main}. Throughout this section, $\Omega$ is the smooth positive closed $(2,2)$-form on $X$ that we define now.

We first observe that the group $\PGL(3,\C)$ acts transitively and holomorphically on $\P^2$. 
Its action lifts to a transitive holomorphic action on $X$. Moreover, $\PGL(3,\C)$ preserves the family of projective lines in $\P^2$.
The projective lines in $\P^2$ lift to disjoint rational curves in $X$ which constitute a smooth holomorphic fibration of $X$. This fibration is invariant under the action of $\PGL(3,\C)$.

Let $L$ be the projective line in $\P^2$ of equation $z=\frac{1}{20} c+ 1$. Let $\Theta$ be a smooth positive form of maximal degree on $\PGL(3,\C)$ supported by a small enough neighbourhood of the identity and of total mass equal to 1. 
Consider the closed positive $(2,2)$-current $\Omega$ on $X$ defined by
$$\Omega:=2\int_{\tau\in \PGL(3,\C)} \tau_*[\widehat L] \Theta(\tau),$$
where $[\widehat L]$ is the current of integration on the lift $\widehat L$ of $L$ to $X$. 

With the above description of the action of $\PGL(3,\C)$ on $X$, it is not difficult to see that $\Omega$ is actually a smooth form. The following two lemmas  are crucial for us. We use here the affine coordinates $(c,z,t)$ for $X$. Then the hypersurface $V\cap \Pi^{-1}(\C^2)$ is given by $t=0$ and we can identify it with $\C^2$ so that $\mu$ is considered as a probability measure with compact support on it. 

\begin{lemma} \label{l:special-form1}
Let $\Kc:=\{(c,z)\in \C^2\;\colon\; c\in\Mc,\ g_c(z)=0\}$. Let  $A$ and $B$ be positive numbers such that $\Kc \subset 
\D_{A}\times \D_{B}$.
For $\delta >0$, let  $T_{\delta,n}$ denote the restriction of ${1\over n2^n} (\widehat F^n)^*(\Omega)$ to 
$\D_{A}\times \D_{B}\times \D_\delta$. Then any limit $T_\delta$ of  $T_{\delta,n}$ is supported on $\Pi^{-1}(\Kc)$.
\end{lemma}

\proof
Let $(c_0,z_0) \in \D_{A}\times \D_{B}\setminus \Kc$. We have to show that there exists a neighbourhood $W$ of 
$(c_0,z_0)$ in $\C^2$ such that $T_\delta$ has no mass on $\Pi^{-1}(W)$. To this end, we will use a form 
$$\Omega_0=2\int_{\tau\in \PGL(3,\C)} \tau_*[\widehat L_0] \Theta_0(\tau)$$ 
which is similar to the form $\Omega$, where $L_0$
is a suitably chosen projective line in $\P^2$ and $\Theta_0$ has a sufficiently small support. We will show that the mass
of ${1\over n2^n} (\widehat F^n)^*(\Omega_0)$ on $\Pi^{-1}(W)$ tends to zero when $n$ tends to infinity. By 
Lemmas \ref{l:coh-only}, \ref{l:Omega-0} and the fact that the dimension of $H^{2,2}(X,\C)$ is $2$, this implies that $T_\delta$ has no mass on $\Pi^{-1}(W)$. We distinguish two cases.

\medskip\noindent
{\bf Case 1.}  Assume that $c_0\notin \Mc$. We have $g_{c_0}(0)>0$. Choose a point $z_0 \in K_{c_0}$. Then we have  $g_{c_0}(z_0)=0$.
Consider $W:=\D(c_0,r_0)\times \D(z_0,\rho_0)$ with $r_0$ and $\rho_0$ small enough so that $g_c(0)\ge 2m_0$ on $\D(c_0,r_0)$ 
and $g_c(z) \le m_0$ on $W$ for some constant $m_0>0$. Observe that $g_c(p_c^n(0))=2^n g_c(0) \ge 2^{n+1} m_0>g_c(z)$ for all $c\in \D(c_0,r_0)$ and $n\ge 1$. Therefore, one sees that $W$ does not meet
the post-critical set of $F$ which is the forward orbit of the line $\{z=0\}$ under the action of $F$.

Choose now $L_0$ and $\Theta_0$ so that $\tau(L_0) \cap (\D(c_0,r_0)\times \C)$ is contained in $W$
for all $\tau$ in the support of $\Theta_0$.
Since for almost every $\tau\in \PGL(3,\C)$  the curve $\widehat{\tau(L_0)}$ does not meet the indeterminacy set of $\widehat F^n $ we have 
\begin{eqnarray}
 (\widehat F^n)^* (\Omega_0)=2\int_{\tau\in \PGL(3,\C)} [\widehat{F^{-n} (\tau(L_0))}] \Theta_0(\tau).
\end {eqnarray}

Then $\tau(L_0)_{\vert_W}$ is the graph  $\Gamma_{\gamma_{\tau}}$ of an affine function  $\gamma_\tau$ over $\D(c_0,r_0)$ which  does not meet the post-critical set of $F$ when $\tau\in \supp\;(\Theta_0)$. It follows that $F^{-n} (\Gamma_{\gamma_\tau})$ is an union of $2^n$ disjoint graphs over $\D(c_0,r_0)$. More precisely, we have 
$$F^{-n} (\Gamma_{\gamma_\tau}) = \bigcup_{j=1}^{2^n} \Gamma_{\gamma_{\tau}^{j,n}} \quad \text{and\ \  hence} \quad 
\widehat{F^{-n} (\Gamma_{\gamma_\tau})} = \bigcup_{j=1}^{2^n} \widehat {\Gamma_{\gamma_{\tau}^{j,n}}},$$
where  $\gamma_{\tau}^{j,n} : \D(c_0,r_0) \to \C$ are holomorphic functions satisfying 
$p_c^n(\gamma_{\tau}^{j,n}(c))= \gamma_{\tau}(c)$ for every integer $n\in \N$ and every $\tau \in \supp\;(\Theta_0)$.

Recall that the function $g_c(z)$ is continuous on $\C^2$ and $\lim_{z\to \infty} g_c(z)=+\infty$ locally uniformly in $c$ (see for instance Lemma \ref{l:Green}(2)). Then, as $g_c(\gamma_{\tau}^{j,n}(c))=\frac{1}{2^n}g_c(p_c^n(\gamma_{\tau}^{j,n}(c)))=\frac{1}{2^n}g_c(\gamma_{\tau}(c))$, one sees that the family
$$\big\{\gamma_{\tau}^{j,n}\;\colon \; \tau \in \supp(\Theta_0),\;  1\le j\le 2^n, n\ge 1\big\}$$ 
is locally uniformly bounded. Thus, after shrinking $r_0$, we have 
$\vert (\gamma_{\tau}^{j,n})' \vert \le M$ on $\D(c_0,r_0)$ for some constant $M>0$ and for all elements of the above family.
Therefore, the graphs $ \Gamma_{\gamma_{\tau}^{j,n}}$ and their lifts $\widehat \Gamma_{\gamma_{\tau}^{j,n}}$ have bounded areas.
This  yields $\Vert {1\over n2^n} (\widehat F^n)^*(\Omega_0)\Vert_{\Pi^{-1}(W)} =O(\frac{1}{n})$.
 
 \medskip\noindent
{\bf Case 2.}  Assume that $z_0\notin K_{c_0}$. Then $g_{c_0}(z_0) >0$ and we may choose a small neighbourhood $W:=\D(c_0,r_0)\times \D(z_0,\rho_0)$ of $(c_0,z_0)$ in $\C^2$  such  that $g_c(z)\ge 2m_0$ on $W$ for some constant $m_0>0$.
We then take $L_0$ and $\Theta_0$ so that $g_c(z) \le m_0$ on $\tau(L_0)\cap (\D(c_0,r_0)\times\C)$
for all $\tau$ in the support of $\Theta_0$.
We have $g_c(z) \le 2^{-n} m_0$ on $F^{-n} (\tau(L_0)) \cap (\D(c_0,r_0)\times\C)$. It follows that the last set is disjoint from $W$ and hence  
$\widehat{F^{-n} (\tau(L_0))} \cap \Pi^{-1}(W)=\varnothing$ for every $\tau\in \supp(\Theta_0)$ and $n\ge 1$.
Hence $\Vert {1\over n2^n} (\widehat F^n)^*(\Omega_0)\Vert_{\Pi^{-1}(W)}=0$ for all $n\ge 1$ and this completes the proof of the lemma.
\endproof

\begin{lemma} \label{l:special-form}
For positive numbers $A$, $B$ and $\delta$, let 
 $T_{\delta,n}$ denote the restriction of ${1\over n2^n} (\widehat F^n)^*(\Omega)$ to 
$\D_{A}\times \D_{B}\times \D_\delta$. If $A, B$ are sufficiently big and $\delta$ is sufficiently small, then 
the support of $T_{\delta, n}$ is contained in $\D_{\frac{A}{2}}\times \D_{\frac{B}{2}}\times \D_\delta$ for every $n\ge 1$
and  $T_{\delta,n}\wedge [V]$ converges to $\mu$ as $n$ tends to infinity.
\end{lemma}

\proof
To start we take any $\delta >0$ and  pick $A$ and $B$ sufficiently big so that $\Kc \subset 
\D_{A}\times \D_{B}$ and $A\ge 4C_0$ where $C_0$ is the constant given by Lemma \ref{l:Green}(3).
In the sequel, we will  have to increase $B$ and decrease $\delta$ in order to get the desired property on the support of 
the current $T_{\delta, n}$.

Let us set $K_{A,R}:=\{(c,z)\in \D_A\times \C\;\colon\; g_c(z)\le R\}$ for a constant $R>0$. 
There exists an $R$ such that $\Pi(\supp(\Omega))\cap (\D_A\times \C)\subset K_{A,R}$ and, by choosing $B$ big enough, we have
$K_{A,R}\subset \D_A\times \D_{\frac{B}{2}}$. Then, as $g_c\circ p_c^n =2^n g_c$, we 
have $F^{-n}(K_{A,R})\subset K_{A,2^{-n}R} \subset K_{A,R}$ and thus
\begin{eqnarray}\label{CondSupp00}
\Pi(\supp \; (\widehat F^n)^*(\Omega))\cap (\D_A\times \C) \subset  \D_A\times \D_{\frac{B}{2}}
\end{eqnarray}
which  in particular implies that 
\begin{eqnarray}\label{CondSupp0}
\supp(T_{\delta,n}) \subset \D_A\times \D_{\frac{B}{2}} \times \D_\delta,\;\forall \delta >0,\forall n\ge 1.
\end{eqnarray}

We now aim to show that $\supp(T_{\delta,n}) \subset \D_{\frac{A}{2}}\times \D_{\frac{B}{2}} \times \D_\delta$ provided that $\delta$ is small enough.
If the support of $\Theta$ is sufficiently small then every $\tau(L)$ has an equation of the form $z=\alpha_\tau c+\beta_\tau$ where 
$10^{-2}\le \vert \alpha_\tau\vert \le 10^{-1}$ and $\vert \beta_\tau\vert \le 2$ for every $\tau$ in the support of $\Theta$ (recall that $L$ is given by $z={1\over 20} c+ 1$ and hence $\alpha_\tau =\frac{1}{20}$, $\beta_\tau =1$ when $\tau$ is the identity). It then follows from our choice of $A$ and Lemma \ref{l:Green}(3) that 
$$g_c(\alpha_\tau c+\beta_\tau)<g_c(c)=2^{-n}g_c(p_c^n(c))<g_c(p_c^n(c))$$ 
for $\vert c\vert \ge \frac{A}{4}$ and $\tau\in \supp(\Theta)$.

Let us denote by $U_c$ the disc $\D(c,\frac{\vert c\vert}{2})$ and, for every $\tau \in \supp(\Theta)$, let us denote by $\gamma_\tau(u)$ the restriction of the function 
$\alpha_\tau u+\beta_\tau$ to the disc $U_c$.
For $\vert c\vert  \ge \frac{A}{2}$ one has $U_c \subset \{\vert u\vert \ge \frac{A}{4}\}$ and thus, according to  what we have just seen, we have
\begin{eqnarray}\label{estau}
g_u(\gamma_\tau(u))<g_u(u)<g_u(p_u^n(u)) \textrm{ when } u\in U_c \textrm{ and } \vert c\vert  \ge \frac{A}{2} \cdot
\end{eqnarray}

In particular, these inequalities show that $\tau(u)\not=p_u^n(u)$ for every $n\geq 0$ and hence
$\Gamma_{\gamma_\tau}$ does not meet the postcritical set of $F$. So we may describe $F^{-n}(\Gamma_{\gamma_{\tau}})$ and its lift
$\widehat{F^{-n}(\Gamma_{\gamma_{\tau}})}$ as
$$F^{-n}(\Gamma_{\gamma_{\tau}})=\bigcup_{j=1}^{2^n} \Gamma_{\gamma_{\tau}^{j,n}} \quad \text{and}\quad
\widehat{F^{-n}(\Gamma_{\gamma_{\tau}})}=\bigcup_{j=1}^{2^n} \widehat{\Gamma_{\gamma_{\tau}^{j,n}}}$$
where $\gamma_{\tau}^{j,n} : U_c\to \C$ are holomorphic and $p_u^n(\gamma_{\tau}^{j,n}(u)) =\gamma_{\tau}(u)$ for all $u\in U_c$ and all $n\ge 1$.
Using again (\ref{estau}) and Lemma \ref{l:Green}(1)(2),  we  get the following estimate for $\vert c\vert \ge \frac{A}{2}$,  $u\in U_c$, $n\ge 1$ and $1\le j\le 2^n$:
\begin{eqnarray*}
\ln \frac{\vert \gamma_{\tau}^{j,n}(u)\vert}{4} &\le& \max\Big(g_u(\gamma_{\tau}^{j,n}(u)),\frac{1}{2} g_u(u)\Big)
= \max\Big(\frac{1}{2^n}g_u(\gamma_\tau(u)),\frac{1}{2} g_u(u)\Big)\\
&\le& \max\Big(\frac{1}{2^n}g_u(u),\frac{1}{2} g_u(u)\Big)
\le \frac{1}{2}\big( \ln 2 +\ln \vert u\vert\big)
\le \frac{1}{2} \ln (3\vert c\vert).
\end{eqnarray*}

Applying  Cauchy inequality on the disc $U_c$, the above estimate  yields some constant $M>0$ such that $\vert (\gamma_\tau^{j,n})'(c)\vert \le M$ for $n\ge 1$, $1\le j\le 2^n$ and $\vert c\vert \ge \frac{A}{2}$. This implies  that   if $N_V$ is a sufficiently thin
neighbourhood of $V$ then $\widehat{F^{-n}(\Gamma_{\gamma_{\tau}}})\cap N_V=\varnothing$ for every $\tau \in \supp(\Theta)$ and every $n\ge 1$. As $\Gamma_{\gamma_{\tau}}=\tau(L)\cap (U_c\times \C)$ for $\vert c\vert \ge \frac{A}{2}$, this means
that 
\begin{eqnarray}\label{CondSupp}
\supp\;((\widehat F^n)^*(\Omega)) \cap N_V\cap\Pi^{-1}\Big(\Big\{\vert c\vert \ge \frac{A}{2} \Big\}\times \C\Big)=\varnothing,\;\forall n\ge1.
\end{eqnarray}
Taking (\ref{CondSupp0}) into account, (\ref{CondSupp}) shows that, for $\delta$ small enough,  the desired inclusion $\supp\;(T_{\delta,n}) \subset \D_{\frac{A}{2}}\times \D_{\frac{B}{2}} \times \D_\delta$ is verified for every $n\ge 1$.

We will now compute the limit of $T_{\delta,n}\wedge [V]$ when $n$ tends to infinity.
By the definition of $T_{\delta,n}$ and  (\ref{CondSupp00}), we have 
$$T_{\delta,n}\wedge [V]= \frac{1}{2^n} (\widehat F^n)^*(\Omega)_{\vert_{\Pi^{-1}(\D_A\times\D_B)}} \wedge [V]= \frac{1}{2^n} (\widehat F^n)^*(\Omega)_{\vert_{\Pi^{-1}(\D_A\times\C)}} \wedge [V]$$ 
which by (\ref{CondSupp}) gives
\begin{eqnarray}\label{FO}
T_{\delta,n}\wedge [V]
=\frac{1}{2^n} (\widehat F^n)^*(\Omega)_{\vert_{\Pi^{-1}(\C\times\C)}} \wedge [V].
\end{eqnarray}
Going back to the definition of $\Omega$ this yields
\begin{eqnarray*}
T_{\delta,n}\wedge [V]&=&\frac{2}{n2^n}\int_{\tau\in \PGL(3,\C)} \Big([\widehat{F^{-n}(\tau(L))}]_{\vert_{\Pi^{-1}(\C\times\C)}} \wedge [V]\Big) \Theta(\tau)\\
&=&
\int_{\tau\in \PGL(3,\C)} \Big(\frac{2}{n2^n} \sum_{(c,z)\in S_{n,\tau}} \delta_{c,z}\Big) \Theta(\tau),
\end{eqnarray*}
where $S_{n,\tau}$ is the set of points where $F^{-n}(\tau(L)_{\vert_{\Pi^{-1}(\C\times\C)}})$ has a vertical tangency. See also Lemmas \ref{l:contact} and \ref{l:contact-lift}.

Let us recall that  $\tau(L)_{\vert_{\Pi^{-1}(\C\times\C)}}=\{(c,z)\in \C\times\C\;\colon z=a_\tau(c)\}$ where $a_\tau(c)=\alpha_\tau c+\beta_\tau$ is a degree one polynomial on $\C$.
Thus, $(c,z)\in \C\times \C$ is a point  where
 $F^{-n}(\tau(L)_{\vert_{\Pi^{-1}(\C\times\C)}})$ has a vertical tangency if and only if $(p_c^n)'(z)=0$ and $p_c^n(z)=a_\tau(c)$.
These conditions are equivalent to $p_c^j(z)=0$ for some $0\le j\le n-1$ and  $p_c^n(z)=a_\tau(c)$. So we may rewrite  $S_{n,\tau}$ as 
\begin{eqnarray*}
S_{n,\tau}=\big\{(c,z)\in \C\times \C\;\colon\; p_c^{n-j}(0)=a_\tau(c),\; p_c^j(z)=0 \text{ for some } 0\le j\le n-1 \big\}.
\end{eqnarray*}

It is enough to show that, for a fixed $a_\tau$, one has $\lim_{n\to\infty} \tilde \mu_n =\mu$ where $(\tilde\mu_n)_n$ is the sequence of discrete measures on $\C\times \C$ defined by 
\begin{eqnarray*}
\tilde\mu_n :=\frac{2}{n2^n} \sum_{j=0}^{n-1} \sum_{c\in I_j^n} \sum_{p_c^j(z)=0} \delta_{c,z} \textrm{ \ and \  }  I_j^n:=\big\{c\in \C\;\colon\; p_c^{n-j}(0)=a_\tau(c) \big\}.
\end{eqnarray*}
Note that it follows from (\ref{FO}) that the above measures $\tilde \mu_n$ are all supported in $\D_A\times\D_B$. Recall again that $V\cap\pi^{-1}(\C\times\C)$ is identified to $\C^2$ and $\mu$ can be seen as a probability measure on $V$.

Denote by  $\pi:\C\times\C\to\C$ the projection $\pi(c,z)=c$ and  set 
$$\mu_n:=\pi_* (\tilde\mu_n)=\frac{2}{n2^n} \sum_{j=0}^{n-1} \sum_{c\in I_j^n}  2^j\delta_c.$$
A potential of $\mu_n$ is given by 
$$\frac{2}{n2^n} \sum_{j=0}^{n-1} 2^{j} \ln \vert  p_c^{n-j}(0)-a_\tau(c)\vert =
\frac{2}{n} \sum_{j=0}^{n-1} \frac{1}{2^{n-j}} \ln \vert  p_c^{n-j}(0)-a_\tau(c)\vert = \frac{2}{n} \sum_{k=1}^{n} \varphi_k(c),$$
where $\varphi_k(c):=\frac{1}{2^k}  \ln \vert  p_c^{k}(0)-a_\tau(c)\vert.$
By Lemma \ref{l:Mandel-Green}, $(\varphi_k)_k$ converges to $g_c(0)$ in $L^1_{\loc}(\C)$ as $k$ tends to infinity. Thus,
$$\lim_{n\to\infty} \mu_n =\lim_{n\to\infty} dd^c \Big(\frac{2}{n}\sum_{k=1}^{n} \varphi_k\Big)= 2dd^c g_c(0)=dd^c g_c(c)=m.$$ 
We used here the identity $g_c(c)=g_c(p_c(0))=2g_c(0)$. 
So if $\tilde \mu'$ is any weak limit of $(\tilde\mu_n)_n$, then we have $\pi_*(\tilde\mu')=m$. Therefore,
we can write  $\tilde\mu'=\int \mu'_c\;dm(c)$
for some probability measures $\mu'_c$ on $\{c\}\times \D_B$. To end the proof of the lemma, it thus suffices to show that 
$\mu'_c=\mu_c$ for $m$-almost every $c$.

For this purpose, we will first compute $F^{N*}(\tilde\mu_{n-N})$ for a fixed $N$ smaller than $n$. Note  that $I_{k-N}^{n-N}=I_k^n$ for $N\le k\le n$. We have
\begin{eqnarray*}
F^{N*} (\tilde\mu_{n-N}) &=&\frac{2}{(n-N)2^{n-N}}  \sum_{j=0}^{n-N-1} \sum_{c\in I_j^{n-N}} \sum_{p_c^{j+N}(z)=0} \delta_{c,z}\\
&=&\frac{2}{(n-N)2^{n-N}}  \sum_{k=N}^{n-1} \sum_{c\in I_{k-N}^{n-N}} \sum_{p_c^{k}(z)=0} \delta_{c,z}\\
&=&\frac{2}{(n-N)2^{n-N}}  \sum_{k=N}^{n-1} \sum_{c\in I_{k}^{n}} \sum_{p_c^{k}(z)=0} \delta_{c,z}\\
&=&\frac{2}{(n-N)2^{n-N}}  \sum_{k=0}^{n-1} \sum_{c\in I_{k}^{n}} \sum_{p_c^{k}(z)=0} \delta_{c,z} - \frac{2}{(n-N)2^{n-N}}  \sum_{k=0}^{N-1} \sum_{c\in I_{k}^{n}} \sum_{p_c^{k}(z)=0} \delta_{c,z}.
\end{eqnarray*}

Observe that the number of terms in the last triple sum is 
$\sum_{k=0}^{N-1} 2^k \# I_k^n =N2^{n-1}$. It follows that (when $N$ is fixed and $n$ tends to infinity)
$$2^{-N}F^{N*} (\tilde\mu_{n-N})=\tilde\mu_n+o\Big({1\over n}\Big).$$
Therefore, $2^{-N} F^{N*}(\tilde\mu_{n-N})$ and $\tilde\mu_n$ have the same limits when $n$ tends to infinity.

Consider any limits $\tilde\mu=\int  \mu'_c\;dm(c)$ of  $(\tilde\mu_n)_n$ and $\tilde\mu_N=\int  \mu'_{c,N}\;dm(c)$ of $(\tilde\mu_{n-N})_n$ using a same subsequence of indices. We have $2^{-N} F^{N*}(\tilde\mu_N)=\tilde\mu$ which implies that
$$\int 2^{-N} p_c^{N*} (\mu'_{c,N})\;dm(c)=\int  \mu'_c\;dm(c).$$ 
We deduce that $2^{-N} p_c^{N*} (\mu'_{c,N})=\mu'_c$ for all $N\ge 1$ and $m$-almost every $c$.
Theorem \ref{t:DS} below, applied to $f:=p_c$, implies that $\mu'_c=\mu_c$ for $m$-almost every $c$ and ends the proof of the lemma. In order to apply Theorem \ref{t:DS}, note that the exceptional set of $p_c$ reduces to $\{\infty\}$ when $c\in b\Mc=\supp (m)$ and the measures $\mu'_{c,N}$ are all supported by $\D_B\Subset \C$. 
\endproof

\begin{theorem}[\cite{DSdens}] \label{t:DS}
Let $f$ be a holomorphic endomorphism of $\P^k$ of algebraic degree $d\geq 2$ and equilibrium measure $\mu$. Let $\Ec$ be the exceptional set of $f$ (the maximal totally invariant proper analytic set of $\P^k$). Let $W$ be a neighbourhood of $\Ec$. Then for all probability measures $\mu_n$ with support in $\P^k\setminus W$ we have 
$$\lim_{n\to\infty} d^{-kn} (f^n)^*(\mu_n)=\mu.$$
In particular, if $\mu'$ is a probability measure such that for every $n$ we have $\mu'=d^{-kn}(f^n)^*(\mu'_n)$ for some $\mu'_n$ supported by $\P^k\setminus W$, then $\mu'=\mu$.
\end{theorem}

\proof[{\bf End of the proof of Theorem \ref{t:main}.}]
Consider a sequence $(n_k)$ tending to infinity such that ${1\over n_k2^{n_k}} (\widehat F^{n_k})^*(\Omega)$ converges to some positive closed current $T$. 
Let  $T_{\delta,n_k}$ be the restriction of ${1\over n_k2^{n_k}} (\widehat F^{n_k})^*(\Omega)$  to $\D_A\times\D_B\times \D_\delta$ where $A$, $B$ and $\delta$ are positive numbers so that the conclusion of Lemma \ref{l:special-form}
holds.
By Proposition \ref{p:limit-22}, there is a probability measure $\nu$ on $\P^2$ such that $T=\Pi^*(\nu)$. Denote by $\nu'$ the restriction of $\nu$ to $\D_A\times\D_B$.  It is enough to show that $\nu'=\mu$ because this identity will imply that $\nu'$ has mass 1 and hence $\nu=\nu'=\mu$. 
Recall that we are using the coordinates $(c,z,t)$ and we can consider $\mu, \nu'$ as measures on $V\cap \Pi^{-1}(\C^2)$.

Denote by $T_\delta$ the restriction of $T$, or equivalently of $\Pi^*(\nu')$, to $\D_{A}\times \D_{B}\times \D_\delta$. So $T_\delta$ is the limit of 
$T_{\delta,n_k}$ and, by Lemma \ref{l:special-form1}, $T_\delta$ is supported by $\Pi^{-1}(\Kc)$.
By Lemma \ref{l:special-form}, all currents $T_{\delta,n_k}$  are vertical-like and supported by $\Pi^{-1}(\overline\D_{\frac{A}{2}}\times\overline \D_{\frac{B}{2}})$. By Lemmas \ref{l:special-form}, \ref{l:slice-order} and Remark \ref{r:vertical}, we have
$$\nu'=T_\delta\wedge [V] \triangleright \lim_{k\to\infty} \big(T_{\delta, n_k} \wedge [V]\big) =\mu.$$

Finally, the first identity in the last line shows that $\nu'$ is supported  by $\Kc$. This, together with the relation $\nu'\triangleright \mu$ and Lemma \ref{l:order-C2}(3), imply that $\nu'=\mu$. This ends the proof of the theorem.
\endproof

\end{document}